\documentclass[a4paper,10pt]{article}
\usepackage{amsxtra}
\usepackage{amssymb}
\usepackage{amscd}
\usepackage{fancyhdr}
\usepackage{graphicx}
\newcounter{ttt}
\newcounter{aaa}
\newcommand{\nc}{\newcommand}

\newtheorem{thm}{Theorem}
\newtheorem{pro}{Proposition}
\newtheorem{lem}{Lemma}
\newtheorem{rem}{Remark}
\newtheorem{deff}{Definition}
\newtheorem{cor}{Corollary}
\newtheorem{ins}{Example}
\newenvironment
{en1}
{
\begin{list}
{\bf{(\roman{ttt} )}}
{
\usecounter{ttt}    \leftmargin=15mm    \rightmargin=5mm
\topsep=0cm         \partopsep=0cm      \labelsep=3mm
\itemsep=1mm        \parsep=0mm         \listparindent=0mm
\itemindent=0mm     \noindent           \labelwidth=1cm
}
}
{\end{list}}
\newenvironment
{en3}
{
\begin{list}
{\bf{(\Roman{ttt} )}}
{
\usecounter{ttt}    \leftmargin=15mm    \rightmargin=5mm
\topsep=0cm         \partopsep=0cm      \labelsep=3mm
\itemsep=1mm        \parsep=0mm         \listparindent=0mm
\itemindent=0mm     \noindent           \labelwidth=1cm
}
}
{\end{list}}
\nc{\bea}{\begin{eqnarray}}
\nc{\eea}{\end{eqnarray}}
\nc{\br}{\bar}
\nc{\fr}{\frac}
\nc{\lm}{\lambda}   \nc{\al}{\alpha}    \nc{\bt}{\beta}
\nc{\dt}{\delta}    \nc{\ep}{\epsilon}  \nc{\zt}{\zeta}
\nc{\et}{\eta}      \nc{\thh}{\theta}   \nc{\kk}{\kappa}
\nc{\sg}{\sigma}    \nc{\pr}{\partial}  \nc{\up}{\upsilon}
\nc{\dg}{\dagger}   \nc{\ve}{\varepsilon}
\nc{\gm}{\gamma}
\nc{\gmm}{\raisebox{.4ex}{$\! \gamma \!$}}
\nc{\st}{\star}
\nc{\stt}{\star \star}
\nc{\mi}{\! \! \mid}
\nc{\vs}{\vskip 1mm}
\nc{\vv}{\vskip .1mm}
\nc{\hs}{\hspace{2mm}}
\nc{\hhs}{\hspace{2em}}
\nc{\acc}{\\[0.5mm]}
\nc{\ul}{\underline}
\nc{\ol}{\overline}
\nc{\cl}{\centerline}
\nc{\lb}{\linebreak}
\nc{\pa}{\parallel}             \nc{\pe}{\perp \!}
\nc{\wt}{\widetilde}                \nc{\df}{\mapsto}
\nc{\td}{\tilde}
\nc{\la}{\leftarrow}                \nc{\Ra}{\Rightarrow}
\nc{\La}{\Leftarrow}                \nc{\ra}{\rightarrow}
\nc{\lla}{\longleftarrow}           \nc{\lra}{\longrightarrow}
\nc{\Lla}{\Longleftarrow}           \nc{\Lra}{\Longrightarrow}
\nc{\lLa}{\longLefttarrow}          \nc{\lRa}{\longRightarrow}
\nc{\LLa}{\LongLefttarrow}          \nc{\LRa}{\LongRightarrow}
\nc{\llra}{\longleftrightarrow}
\nc{\Llra}{\Longleftrightarrow}
\nc{\B}[1]{{\bold#1}}
\nc{\bs}[1]{{\boldsymbol#1}}
\nc{\mvb}{\mathversion{bold}}
\nc{\mvn}{\mathversion{normal}}
\nc{\C}[1]{{\mathcal#1}}
\nc{\D}[1]{{\mathbb#1}}
\nc{\F}[1]{{\mathfrak#1}}
\nc{\UU}{\C{U}}
\nc{\DD}{\C{D}}
\nc{\OO}{\C{O}}
\nc{\cc}{\D{C}}
\nc{\hh}{\D{H}}
\nc{\nn}{\D{N}}
\nc{\oo}{\D{O}}
\nc{\qq}{\D{Q}}
\nc{\rr}{\D{R}}
\nc{\zz}{\D{Z}}
\nc{\xx}{\F{X}}
\nc{\sst}{\scriptstyle}
\nc{\sss}{\scriptscriptstyle}
\nc{\dx}[2]
{
\setbox0=\hbox{\raise.1em\hbox{#1}\kern-1.1em\lower.4em\hbox{#2}}
\box0 \!
}
\nc{\wa}{\! \! \!}
\nc{\wb}{\wa \wa}
\nc{\wc}{\wa \wb }
\nc{\zzz}{\wa}
\nc{\ci}{\hfill{\bf \bigodot}}
\nc{\fpum}{\Phi _r(\perp ^{1}\wb M)}
\nc{\fpus}{\Phi _r(\perp ^{1}\wa S)}
\nc{\pem}{\perp \wa  M}
\nc{\pes}{\perp \wa  S}
\nc{\pus}{\perp \! ^{\wa 1}S}
\nc{\pum}{\perp \! ^{\wa 1}M}
\nc{\pun}{\perp \! ^{\wa 1}N}
\nc{\ip}[2]{\langle #1,#2 \rangle }
\nc{\btab}{\begin{tabular}}
\nc{\etab}{\end{tabular}}
\nc{\bite}{\begin{itemize}}
\nc{\eite}{\end{itemize}}
\nc{\ben}{\begin{en1}}
\nc{\een}{\end{en1}}
\nc{\barr}{\begin{array}}
\nc{\earr}{\end{array}}
\nc{\bca}{\begin{cases}}
\nc{\eca}{\end{cases}}
\nc{\bthm}{\begin{thm}}
\nc{\ethm}{\end{thm}}
\nc{\bpro}{\begin{pro}}
\nc{\epro}{\end{pro}}
\nc{\blem}{\begin{lem}}
\nc{\elem}{\end{lem}}
\nc{\bins}{\begin{ins}}
\nc{\eins}{\end{ins}}
\nc{\bcor}{\begin{cor}}
\nc{\ecor}{\end{cor}}
\nc{\brem}{\begin{rem}}
\nc{\erem}{\end{rem}}
\nc{\bdeff}{\begin{deff}}
\nc{\edeff}{\end{deff}}
\def\tdnabla{\wt \nabla}

\def\pf{{\bf{Proof:} \\ }}
\def\bpf{\pf}
\def\epf{$\ci$}
\def\go{{\bf{($\Lra$)} \\ }}
\def\come{{\bf{($\Lla$)} \\ }}

\begin{document}

\title{Superminimal Surfaces in the 6-Sphere}
\author{Martins, J.K.\footnote{Work fully supported by CAPES, FAPEAM and Universidade Federal do Amazonas.}}
\maketitle

\begin{abstract}

In this article, we use the harmonic sequence associated to a
weakly conformal harmonic map $f:S\to S^6$ in order to determine
explicit examples of linearly full almost complex 2-spheres of
$S^6$ with at most two singularities.  We prove that the
singularity type of these almost complex 2-spheres has an extra
symmetry and this allows us to determine the moduli space of such
curves with suitably small area. We also characterize
projectively equivalent almost complex curves of $S^6$ in terms
of $G_2^{\cc}$-equivalence of their directrix curves.

{\bf Email:} akay333@gmail.com

\end{abstract}

\section{Preliminaries}
The use of harmonic sequence is a very well known technique since,
in recent times, it has been used by several authors
(\cite{dz1},\cite{ew1},\cite{bw1}) and we shall give here the same
treatment as in \cite{bw1}, which was suitably specialized for the
case of the spheres. We start by establishing some background
results about almost complex curves of the 6-sphere.

A map $\psi : S\to W$ between Riemannian manifolds is {\bf
harmonic} if it satisfies the {\bf Euler-Lagrange equation}
$$
{\rm tr}(\nabla d\psi )= 0.
$$
Throughout this article, we shall use $S$ to denote a Riemann
surface and $z=x+iy$ shall denote a local complex coordinate z on
$S$. In this case, the harmonicity condition of $\psi$ is
simplified to

\begin{equation} \label{su1}
(\nabla_{\frac{\partial}{\partial\bar{z}}}d\psi)
(\frac{\partial}{\partial z}) = 0.
\end{equation}

Let $V \to S$ be a complex vector bundle over the Riemann surface
$S$ and assume that $\tdnabla$ is a connection on $V$. By the
Koszul-Malgrange Theorem, $V$ admits the structure of a
holomorphic vector bundle. Here a section $s$ is a {\bf
holomorphic section} if and only if
\begin{equation} \label{su2}
\tdnabla_{\frac{\partial}{\partial\bar{z}}}s=0.
\end{equation}

Given a harmonic map $\psi _0:=\psi : S \to \cc P^n$, several
authors (\cite{wolfson1},\cite{ew1},\cite{dz1},\cite{bw1}) have
dealt with the sequence of harmonic maps $\psi _k: S \to \cc P^n$
obtained from $\psi$ via an inductive construction of a sequence
of complex line bundles over $S$.

In the sequel, we outline this construction and give some of the
main features of this sequence.

Let $\F L$ be the tautological line bundle over $\cc P^n$. Let
$L_0$ and $L_0^\perp$ be the pullbacks via $\psi _0$ of $\F L$ and
$\F L^\perp$ respectively. $L_0$ and $L_0^\perp$ are endowed with
naturally induced connections for they are vector subbundles of
the trivial $\cc^{n+1}$-bundle over $S$.

Explicitly, if $s$ is a section of a subbundle $L$ of the trivial
bundle $S\times \cc^{n+1}$, then $s$ may be regarded as a map $S
\to \cc^{n+1}$. Given $X \in T_xS$, we can define a connection
$\nabla_Xs$ by the orthogonal projection $(Xs)^{L}$ of $Xs$ onto
$L$. Similarly, we also define a connection in $L^\perp$. Thus the
line bundles $L_0$ and $L_0^\perp$ have structures of holomorphic
vector bundles over S.

The map $\psi _0$ determines a bundle map $\partial_0 : L_0 \to
L_0^\perp$. Indeed, if we consider a holomorphic local section
$f_0 : S \to \cc^{n+1}\setminus\{0\}$ of $L_0$, then we define
$\partial_0f_0 = (\frac{\partial f_0}{\partial z})^{L_0^\perp},$
and $\overline{\partial}_0 f_0 = (\frac{\partial f_0}{\partial
\bar{z}})^{L_0^\perp}$.
%
%

It follows from \eqref{su1} and \eqref{su2} that $\partial_0$
($\ol{\partial}_0$) is a holomorphic (anti-holomorphic) bundle map
if and only if $\psi_0 = [f_0]$ is a harmonic map. Therefore if
$\td f_1:=
\partial_0f_0$ is not identically zero then it is a holomorphic section
of the bundle $L_0^\perp$ and hence its zeros (if any) are
isolated. Let $z_0$ be such a zero, then for some holomorphic
local section $f_1$ we have $\td f_1(z) = (z-z_0)^rf_1(z)$ with
$f_1(z_0) \ne 0$. This latter map will then yield a well-defined
map $\psi_1(z) := [f_1(z)]$ from $S$ into $\cc P^n$ and $f_1$ is a
meromorphic local section for a complex line bundle $L_1\subset
L_0^\pe$.

By defining $\partial_1$ in a similar way, and verifying that $\td
f_2 := \partial_1 f_1$ is a holomorphic local section of
$L_1^\perp$, it follows that $\psi _1$ is also harmonic.

Therefore, as long as the bundle section $f_k$ is not identically
zero, that is, $\psi _{k-1}$ is not anti-holomorphic (or $\psi
_{k+1}$ is not holomorphic, when considering the descending
sequence given by $\psi _{-i-1}=\ol{\partial}_{-i}(\psi _{-1})$
where $i=0,1,2,\hdots $), we can carry on with this
process, defining a sequence $\psi_k = [f_k(z)]$ of harmonic maps
such that the local sections $f_k$ are characterized by the
following properties:

\bea \frac{\partial f_p}{\partial z} &=&
f_{p+1} + \frac{1}{|f_p|^2}\ip{\frac{\partial}{\partial z}f_p}{f_p}f_p
\notag \\
&=& f_{p+1} + \partial _z(log|f_p|^2)f_p
\notag \\
&=& f_{p+1} + \al _pf_p \mbox{\hs \hs where \hs \hs } \al _p:=
\partial _z(log|f_p|^2)
\label{su3} \\
\frac{\partial f_{p+1}}{\partial \overline{z}} &=&-\gm _pf_p
\phantom{\al _pf_px} \mbox{\hs \hs where \hs \hs } \gm _p:=
\frac{|f_{p+1}|^2}{|f_p|^2}. \label{su4} \eea

Here $\ip{}{}$ denotes the standard Hermitian product on $C^{n+1}$
and $|\hs |$ denotes the associated norm. It is known
\cite{ew1} that the harmonic sequence $\psi _i$ where $i\in\zz$,
terminates at one end if and
only if it terminates at both ends. If this happens, we say that
each element of the sequence is \ul{\bf superminimal} and it is
customary to consider the range for the indices starting at the
holomorphic element of the sequence, that is, $\big(
\psi_j\big)_{j=0}^n$ denotes the harmonic sequence generated by $\psi_0$. This
holomorphic map is usually named in the literature as the \ul{\bf
directrix curve} associated to any of the harmonic maps $\psi_j$ referring to the terminology
adopted when dealing with harmonic 2-spheres in $S^{2m}$(cf.
\cite{barbosa1}).

If $\psi _m=[f_m]$ for some harmonic map $f_m:S \to S^n$, then we
can consider $f_m$ as a nowhere vanishing global holomorphic
section of $L_0$ so that the sequence of meromorphic sections
$f_j$ will also satisfy the condition:

\bea
\ol f_{m+k} &&=
(-1)^k|f_{m+k}|^2f_{m-k}.\label{su5}
\eea

In particular,

\bea
|f_{m+k}| |f_{m-k}|\equiv 1.\label{su5a}
\eea

Note that in this
situation, the element $\psi _0$ will necessarily be in the middle
of the sequence, that is $n=2m$, because $f_{m+k}\equiv 0$ if and
only if $f_{m-k}\equiv 0$.  Moreover, \eqref{su3} and $|f_m|\equiv
1$ implies \bea
\partial _zf_m=f_{m+1}. \hs \hs \hs
\big(Notation\hs
\partial _z
:=\frac{\partial }{\partial z}
 =\frac{1}{2}(\frac{\partial }{\partial x}
  -i\frac{\partial }{\partial y})
\big] \label{su5b} \eea

\bdeff \label{stsu1}
We say that a map
from a Riemannian manifold $N$ into $\cc P^n$ is \ul{\bf linearly
full}, when its image is not contained in any complex space form
$\cc P^k$ for $k<n$.
\edeff

If the Riemann surface is homeomorphic
to the sphere $S^2$ then Wolfson \cite{wolfson1} shows that the
corresponding complex line bundles are mutually orthogonal and
consequently the harmonic sequence terminates, that is, all the
harmonic 2-spheres of $\cc P^n$ are superminimal. Moreover, in
this case, the length of the sequence achieves its maximum, $n+1$,
if and only if $\psi$ is linearly full.

A detailed discussion of the holomorphic curves of a complex
projective space can be found in \cite{gh} (pages 263-268), but
here we describe some material on this topic to be used in this
article.

Let $\psi (z)=[f(z)]:S \to \cc P^n$ be a holomorphic curve from
the Riemann surface $S$ into $\cc P^n$, where $f:S \to \cc
^{n+1}\setminus\{0\}$ is a local holomorphic lift of $\psi$. Then
the {\bf $\bs{j^{th}}$-osculating curve} of $\psi$ is the
holomorphic curve $\sg _j:S \to \cc P^{n_j}$ (\hs where \hs \hs $
n_j:= \left[\begin{smallmatrix} n+1\\ j+1 \end{smallmatrix}
\right] -1 $) defined by
$$
\sg_j(z)=[f \wedge \hdots \wedge f^{(j)}](z),
$$
where $f^{(j)}=\frac{\pr ^jf}{\pr z^j}$ and $j=0,...,n-1$.

A {\bf higher order singularity} of $\psi$ is a point $p\in S$
which is a singularity for some $j^{th}$-osculating curve
$(j=0,1,..,n-1)$. The {\bf ramification index} of $\sg_j$ at a
point $p$ is the order $r_{j+1}$ of this point as a zero of the
derivative of the curve $\sg_j$.

Thus, the holomorphic curve $\psi$ is said to have {\bf
singularity type} $(r_1,r_2,,..,r_n)$ at the point $p$.

If $S$ is a compact Riemann surface, then the curve $\psi$ has a
finite set $\C Z_{\sss \psi}=\{ p_1,..,p_k\}$ of higher order
singularities. We shall denote by $R_{j+1}$ the sum of the
ramification indices of $\sg _j$ at each singularity, that is,
\bea R_{j+1}=\sum_{i=1}^{k}r_{j+1}(p_i) \label{su34} \eea and we
shall refer to $R_{j+1}$ just as the ramification index of $\sg
_j$.  Moreover, we can define the {\bf total ramification index}
of $\psi$ as the sum $\sum_{j=1}^nR_j$.

Then, we shall say that the holomorphic curve $\psi$ or any
element of its corresponding harmonic sequence has {\bf
singularity type} $(r_1(p),\hdots ,r_n(p))$ at the point $p$ and
has {\bf total singularity type} $(R_1,\hdots ,R_n)$.

The curve $\psi$ is {\bf totally unramified} when its total
ramification index is zero. Otherwise, $\psi$ is said to be {\bf
$\bs k$-point ramified} if the set $\C Z_{\sss \psi}$ has
cardinality $k$.

In terms of a local complex coordinate $z$ for $S$ centred on $p$,
that is $z(p)=0$, it is possible to determine a basis $\{
v_0,..,v_n\}$ for $\cc ^{n+1}$ in such a way that the holomorphic
map $f$ can be written in the normal form
\begin{align}
f(z)& =\sum_{i=0}^{n}z^{k_0+\hdots +k_i}h_i(z)v_i, \label{su32}
\end{align}
where $k_0=0,\hs k_i=r_i(p)+1$ ($j=1,\hdots ,n$) and $h_i(z)$
denotes a holomorphic function satisfying $h_i(0)\not= 0$.

If $\psi :S^2\to \cc P^n$ is 2-point ramified, then we can find a
local complex coordinate for $S^2$ so that the higher order
singularities of $\psi$ (if any) occur at $z=0$ and $z=\infty$.
Indeed, this follows from the fact that the M\" obius group of
conformal transformations of the $2$-sphere acts triply
transitively on $S^2$.

We observe that $\psi$ is a holomorphic map between algebraic
varieties, since $S^2=\cc P^1$, and so $\psi$ is an algebraic map
(c.f. \cite{gh}), that is $f$ is a rational function. Without loss
of generality, we may assume $f$ to be a $\cc ^{n+1}$-valued
polynomial function. \bdeff We say that a harmonic map $\psi : S^2
\to \cc P^n$ has {\bf $\bs{S^1}$-symmetry} if there exists
non-trivial $S^1$-actions on $S^2$ and $\cc P^n$, where the action
on $\cc P^n$ is by holomorphic isometries, such that for all $z\in
S^2$,
$$
\psi (e^{i\thh}z)=e^{i\thh}\psi (z).
$$
\edeff

We are now able to state the main results (\cite{bw2},\cite{bvw2})
to be used in the following sections, which are concerned with the
characterization of the $k$-point ramified harmonic 2-spheres of
$\cc P^n$ for $k\leq 2$.

\bthm \label{stsu4}

Let $\psi : S^2 \to \cc P^n$ be a linearly full harmonic map with $S^1$-symmetry.
Then there exists a holomorphic coordinate $z$ on $S^2$ such that the
directrix curve $\psi _0=[f_0]$ of $\psi$ can be expressed up to
holomorphic isometries of $\cc P^n$ by
$$
f_0(z) = \sum_{p=0}^nz^{k_1+\hdots +k_p}v_p, \label{su6}
$$
where $\{ v_0,\hdots ,v_n\}$ is an orthogonal basis of $\cc
^{n+1}$ and the scalars $k_j$ are positive integers. Furthermore,
$\psi$ is either totally unramified or 2-point ramified. In the
latter case $\psi$ has singularities at $z=0$ and $z=\infty$ with
corresponding singularity type at $z=0$ and $z=\infty$ given
respectively by $(k_1-1,\hdots ,k_n-1)$ and $(k_n-1,\hdots
,k_1-1)$. \ethm

\bthm \label{stsu5}

Let $\psi :S^2 \to \cc P^n$ be
a linearly full harmonic map which is $k$-point ramified for
$k\leq 2$. Let $z$ be a complex coordinate on $S^2$, then
\begin{en1}
\item The higher order singularities of $\psi$ (if any) occur at
$z=0$ and $z=\infty$ if and only if its directrix curve $\psi
_0(z)=[f_0(z)]$ can be written in the form
\bea
f_0(z)=\sum_{p=0}^nz^{k_1+\hdots +k_p}v_p, \label{su7}
\eea
where the scalars $k_j$ are positive integers and the vectors $v_j$
constitutes a basis of $\cc ^{n+1}$. Furthermore, $\psi$ has
$S^1$-symmetry (with fixed points of the $S^1$-action at $z=0$ and
$z=\infty$) if and only if the basis $\{ v_0,\hdots ,v_n\}$ is
orthogonal. \item $\psi (S^2)\subset \rr P^n$ and one of the
conditions (and so both) in the first equivalence stated in {\bf
(i)} occurs if and only if $n=2m$ for some integer $m$ and the
following properties are satisfied.
\bea
k_{j}&=&k_{n-j+1} \mbox{\hspace{18mm} for $i,j\in \{1,\hdots ,n\}$,} \label{su9} \\
\ip{v_j}{\ol{v}_i} &=& (-1)^{j}\dt _{\sss (j,n-i)}\mu\lm _j,
\mbox{\hs \hs for $i,j\in \{0,\hdots ,n\}$,} \label{su9a}
\eea
where $\mu$ is a constant and $ \lm _j:= \frac{\prod_{1\leq r\leq
s\leq n}(k_r+\hdots +k_s)} {\prod_{r=1}^{j}(k_r+\hdots +k_j)
\prod_{r=1}^{n-j}(k_{j+1}+\hdots +k_{n-r+1})}. $
\end{en1}
\ethm

The constant $\mu$ in the theorem can be chosen to be $1$ by
rescaling the homogeneous coordinates of $\cc P^n$. However, we
will avoid this in order to facilitate our calculations later on
when determining examples of superminimal almost complex curves.

\brem \label{stsu5a}
The theorem above shows in particular that
there does not exist a 1-point ramified linearly full harmonic
2-sphere in $\cc P^n$.
\erem

\bdeff \label{stsu5a1}

Two maps $\psi
,\wt \psi:S \to \cc P^n$ are said to be projectively equivalent if
there exists $[A]\in PGL(n+1,\cc )$ so that $\wt \psi =[A](\psi
)$.

\edeff

\bcor \label{stsu5b}

Any two $k$-point ramified ($k\leq
2$) linearly full harmonic maps $\psi ,\wt \psi :S^2\to \cc P^n$
with the same singularity type are projectively equivalent up to
a conformal transformation of $S^2$.

\ecor

\bpf
According to Theorem \eqref{stsu5} these curves are uniquely
determined by their singularity type and a choice of basis for
$\cc ^{n+1}$. Thus, item (i) of that theorem shows that the
corresponding directrix curves differ by an element of $GL(n+1,\cc
)$.\epf

\section{The twistor fibration $\pi :Q^5 \to S^6$}

Let $Q^5$ denote the quadric of $\cc P^6$, which is the K\" ahler
submanifold defined
$$
Q^5=\{ [x]\in \cc P^6 \mbox{\hs such that \hs $(x,x)=0$}\},
$$
where $(,)$ denotes the symmetric bilinear Euclidean inner product on $C^7$ given by
$(x,y) = \langle x, \overline y \rangle$.

The twistor fibration $\pi :Q^5 \to S^6$ is defined by the map
$$
\pi ([x])=\frac{i}{|x|^2}\ol{x}\times x.
$$
Where the cross product $(\cc^7,\times )$ we are considering here
is given in \cite{bvw1} and we shall also be using its main
properties presented in that reference.

We shall see in the next section that the superminimal almost
complex curves of the 6-sphere can be characterized as the
projections of a special type of holomorphic curves of this
quadric. This led us to investigate what is the group of
holomorphic transformations of the quadric which preserves the
superhorizontal distribution to be defined ahead.

Although $\pi$ is not a Riemannian submersion, it is quite close
to that as we shall see below. Moreover, the fact that $\pi$ can
be easily expressed in terms of the cross product $\times$ on $\rr
^7$ (extended $\cc$-linearly to $\cc ^7$) yields some good methods
to investigate properties of any lifting.

The exceptional Lie group $G_2$ acts transitively on the manifolds
$Q^5$ and $S^6$ in such a way that these manifolds can be realized
also as the homogeneous spaces $G_2/U(2)$ and $G_2/SU(3)$
respectively. By considering
these homogeneous spaces, it is possible to show
that the twistor fibration just defined is nothing but the
canonical projection of the first space onto the second one.

We can write any element of $Q^5$ as
$[x]=[a-ib]$ where $a$ and $b$ are orthonormal vectors of $\rr
^7$. In this case $\pi$ reduces to
\begin{equation}\label{3c2}
\pi [x]=a\times b.
\end{equation}

\brem \label{stkm11}

It follows from the characterization of $G_2$
as the group of automorphisms of $(\rr ^7,\times )$ that the map
$\pi$ is $G_2$-equivariant, that is, $\pi[gx]=g(\pi [x])$ for
every $g \in G_2$.

\erem

Using the Hopf fibration $H:S^{13}\to CP^6$, we have for each $x\in S^{13}\subset\cc^7$,
an isomorphism $dH_x:T'_xS^{13}\to T_{[x]}\cc P^6$. Where the horizontal space:

$$
T'_xS^{13}= \{ v\in \cc ^7 \mbox{\hs such that \hs $(v,\ol
x)=0$} \} =\{ix\}^{\perp}_{\cc}
$$
induces the natural decomposition
$T_xS^{13}=span_{\rr}\{ix\}\oplus T'_xS^{13}$
and it yields also the isomorphism:
$$
T_{[x]}Q^5  \cong \{ v\in \cc ^7 \mbox{\hs such that \hs } (v,\ol x)=0
\mbox{\hs and \hs} (v,x)=0. \}
$$

From now on, we shall be identifying these tangent spaces with the linear subspaces of $\cc^7$ mentioned above with no further reference. Using the definition of $\pi$ we have

\begin{equation} \label{3c3}
\pi _{*}|_{[x]}(v)=\frac{i}{|x|^2}.(\ol{x}\times v - x\times
\ol{v}).
\end{equation}

Thus, the vertical distribution on $Q^5$ defined as the kernel of
$\pi _*$ is given by the set of those tangent vectors $v\in
T_{[x]}Q^5$ such that the imaginary part of $x\times \ol{v}$ is
zero. However, writting $v=c+id$, using the notation $\Im$ and
$\Re$ to denote the imaginary vector and real vector parts in $\rr^7$ of a vector in $\cc^7$, we have
\bea \Im (x\times \ol{v})\times d &=(d,d)a-(b\times c)\times d
&= - \Re (x\times \ol{v})\times c \notag \\
\Im (x\times \ol{v})\times c &=(c,c)b-(a\times d)\times c &=  \Re
(x\times \ol{v})\times d. \notag \eea
In order to obtain the equations above, we have used the ordinary property:
$$
u\times (v\times w)+(u\times v)\times w=2(u,w)v-(u,v)w-w(v,u).
$$
Thus, $\Im (x\times
\ol{v})=0$ if and only if $\Re (x\times \ol{v})=0$ which implies
that the vertical distribution is characterized by

\begin{equation} \label{3c4}
V_{[x]}= Ker\pi _*=\{ v\in T_{[x]}Q^5 \mbox{ such that } x \times
\ol{v}=0\}.
\end{equation}

We should also note that $V_{[x]}$ is an isotropic subspace of
$T_{[x]}Q^5$, since we have for any $v\in V_{[x]}$:
$$
0=\ol{v} \times (x\times \ol{v})=2(\ol{v}, \ol{v})x.
$$
This yields a distribution of isotropic subspaces $H':=\ol{V}$ of
the horizontal spaces $H=V^\pe$, which we henceforth will name as
the {\bf superhorizontal distribution}. It follows then from
\eqref{3c4} that this vector space is characterized at $[x]$ by:
\begin{equation} \label{3c5}
H'=\{ v\in T_{[x]}Q^5 \mbox{ such that } x \times v=0 \}.
\end{equation}
Incidentally, looking at the point $v_0=\pi [x] \in S^6$ as a real
vector of $\cc ^7$, it is clear that $v_0\in T_{[x]}Q^5$.
Moreover, $v_0$ is a horizontal vector because we have for any
vertical vector $v\in V$
$$
\ip{\pi [x]}{v}=(\ol{x}\times x,\ol{v})=-(\ol{v}\times
x,\ol{x})=0. \hs \hs \forall \hs v\in V.
$$
Furthermore, the equation above also shows that $v_0$ is
orthogonal to $H'$. Thus, using the 1-dimensional complex space
$D$ spanned by $v_0$, we can split the horizontal distribution as
follows.
\begin{equation} \label{3c6}
H=D\oplus H'.
\end{equation}
We shall investigate now how far the map $\pi$ is prevented from
being a holomorphic Riemannian submersion, by looking at its
behaviour concerning to length-preservation and $\cc$-linearity of
its differential.  We split these properties into two lemmas.

\blem \label{stkm12} $\pi _*$ is length-preserving in $H'$ and it
reduces the length by a $\sqrt{2}$-factor in $D$. \elem

\blem \label{stkm13} $\pi _*$ is $\cc$-linear in $H'$ and it is
$\cc$-anti-linear in $D$. \elem

We have that $G_2$ acts on $Q^5$ as a transitive group of isometries which preserve
horizontal subspaces. Hence it is enough to show both Lemmas at the point
$x=e_1 + ie_5$. In this case, we have $\pi([x])=e_4$ and the vertical subspace is given by
$V_{[x]} = {\rm span}_{\bf C} \{e_2 + ie_6, e_3 + ie_7\}$ and $H'=\ol V_{[x]}$.

{\bf Proof of Lemma \eqref{stkm12}: }
The restriction $d\pi|_{H'}$ is length preserving as for
$v = e_2 - ie_6$ or $v=e_3 - ie_7$ we have $|v| = |x| = \sqrt{2}$
giving $\|v\| = 2$ in the Fubini-Study metric and $d\pi(v) = 2e_7$ or $2e_6$,
giving $|d\pi(v)| = 2 = \|v\|$.

The restriction $d\pi|_{D}$ is reducing length since $|e_4| = 1$, $|x| = \sqrt{2}$
so $\|v\| = \frac{2\cdot 1}{\sqrt{2}} = \sqrt{2}$ in the Fubini-Study metric and $d\pi(e_4) = -e_1$
has length $|d\pi(e_4)| = |e_1| = 1 = \frac{1}{\sqrt{2}}\|e_4\|$.
So $d\pi$ reduces lengths by the factor $\sqrt{2}$. \epf

{\bf Proof of Lemma \eqref{stkm13}: }
The map $d\pi|_{H'}$ is complex linear, since for
$v = e_2 - ie_6$ or $v=e_3 - ie_7$ we have
$d\pi(iv) = Jd\pi(v) = e_4 \times d\pi(v)$.

The map $d\pi|_{D}$ is complex anti-linear, since
$d\pi(ie_4) = -Jd\pi(e_4) = - e_4 \times d\pi(e_4)$. \epf

Hopf hypersurfaces of space forms have been studied in \cite{martins2}
and the twistor fibration presented here, can
provide a way to study this type of hypersurface in the quadric
$Q^5$. The results above show that although $\pi$ is not a
Riemannian submersion it is not so far from this. Consequently, we
can ask whether the lift $\wt M=\pi ^{-1}(M)$ of a Hopf
hypersurface $M$ of $S^6$ is also a Hopf hypersurface in $Q^5$ or
not. However, it is not hard to see that the horizontal lift of a
normal vector field on $M$ cannot lie either in the distribution
$D$ or in $H'$. This fact makes clear that $\wt M$ cannot be a
tubular hypersurface. Furthermore, the decomposition \eqref{3c6}
of the horizontal distribution makes it rather complicated to work
with the Riemannian connection of $Q^5$.

We shall see later on in this article that the superminimal almost
complex curves of $S^6$ are in 1-1 correspondence with the
holomorphic curves of $Q^5$ which are tangential to the
superhorizontal distribution. This motivates us to determine what
is the group of holomorphic transformations of $Q^5$ which
preserves the superhorizontal distribution.

Let us consider the Lie group $ H_1= \{ \lm I \in GL(n+1,\cc ):\lm
\in \cc ^* \} $ and its Lie subgroup $ H_2= \{ \lm I \in
SO(n+1,\cc ): \lm \in \cc \mbox{\hs and \hs} \lm ^{n+1}=1 \} $. It
is well known (for instance, \cite{gh} page 65) that $PGL(n+1,\cc
)=GL(n+1,\cc )/H_1$ is the group of holomorphic transformations of
$\cc P^n$. The following lemma is easy to prove.

\blem

$SO(n+1,\cc )/H_2$ is the group of holomorphic
transformations of the quadric $Q^{n-1}$.

\elem

%
%
%

In the next proposition we will need the following elementary
properties of the distributions $D,H'$ and $V$.
$$
\begin{array}{ll}
D_{[x]}&=D_{[\ol{x}]} \\
H'_{[x]}&=V_{[\ol{x}]} \\
V_{[x]}&=H'_{[\ol{x}]}.
\end{array}
$$

\brem \label{st3c1}

Let $G_2^{\cc}$ be the group of automorphisms
of $(\cc ^7,\times )$. It is clear from the definition of $G_2$
and the characterization of the superhorizontal distribution given
in \eqref{3c5}, that $G_2^{\cc}$ is a Lie subgroup of $SO(7,\cc )$
and it preserves the superhorizontal distribution.

\erem

A {\bf ${\bf G_2}$-basis} (or {\bf canonical basis}) of ${\bf R}^7$ is an orthonormal
basis $\{f_1,\ldots,f_7\}$ of ${\bf R}^7$ satisfying the relations:
$ f_3 = f_1 \times f_2, f_5 = f_1 \times f_4, f_6 = f_2 \times f_4,	f_7 = f_3 \times f_4.$
Hence, if $f_1,f_2,f_4$ are orthogonal unit vectors such that $f_4 \perp f_1 \times f_2$
then $f_1,f_2,f_4$ determine a unique $G_2$-basis subject to the relations
$f_i\times(f_j\times f_k) + (f_i \times f_j) \times f_k = 2\delta_{ik}f_j
	-\delta_{ij}f_k - \delta_{jk}f_i.$
Every $G_2$-basis will have the following multiplication table
\begin{eqnarray*}
		f_i \times f_j \; = \;
		\begin{array}{r||c|c|c|c|c|c|c}
		i\backslash j & 1 & 2 & 3 & 4 & 5 & 6 & 7 \\ \hline \hline
		1 & 0 & f_3 & -f_2 & f_5 & -f_4 & -f_7 & f_6 \\ \hline
		2 & -f_3 & 0 & f_1 & f_6 & f_7 & -f_4 & -f_5 \\ \hline
		3 & f_2 & -f_1 & 0 & f_7 & -f_6 & f_5 & -f_4 \\ \hline
		4 & -f_5 & -f_6 & -f_7 & 0 & f_1 & f_2 & f_3 \\ \hline
		5 & f_4 & -f_7 & f_6 & -f_1 & 0 & -f_3 & f_2 \\ \hline
		6 & f_7 & f_4 & -f_5 & -f_2 & f_3 & 0 & -f_1 \\ \hline
		7 & -f_6 & f_5 & f_4 & -f_3 & -f_2 & f_1 & 0
		\end{array}
\end{eqnarray*}
For any two $G_2$-bases $\{f_1,\ldots,f_7\}$ and $\{\tilde{f}_1,\ldots,\tilde{f}_7\}$
there exists a unique $g \in G_2$ such that $gf_i = \tilde{f}_i$ (we simply define $g$ by
$gf_i = \tilde{f}_i$ and check that it is in $G_2$). Also $g \in G_2$ mapping $e_1,\ldots,e_7$
to $f_1,\ldots,f_7$ can be represented by the matrix
\[
	(f_1|\ldots|f_7) \in SO(7)
\]
(with respect to the standard basis $e_1,\ldots,e_7$). Further details on $G_2$-basis can
be found in \cite{bvw1}.

\bpro
The group $\wt G$ of holomorphic transformations of $Q^5$ which
preserves the superhorizontal distribution is $G_2^{\cc}$.
\epro
\bpf
In accordance with the remark \eqref{st3c1} above, we have
$G_2^{\cc}\subset \wt G$. Let $[g]$ be an arbitrary element of
$\wt G$, that is, $[g]\in SO(7,\cc )/H_2$. Thus $+g$ or $-g$ lies
in $SO(7,\cc )$. By assumption we have for every $[x]\in Q^5$ and
$v\in H'_{[x]}$ that \bea gx\times gv=0. \label{3c9} \eea We shall
first observe that $\ol g$ also lies in $\wt G$. Indeed, given
$v\in H'_{[x]}$ we have $\ol{v}\in V_{[x]}=H'_{\ol{x}}$ and hence
\bea \ol{g}x\times \ol{g}v=\ol{g\ol{x}\times g\ol{v}}=0.
\label{3c10} \eea The superhorizontal subspaces at the points
$[x_0]=[e_1-ie_5]$ and $[x_1]=[e_1-ie_4]$ are
$H'_{[x_0]}=span_{\cc}(e_2+ie_6,e_3+ie_7)$ and
$H'_{[x_1]}=span_{\cc}(e_7+ie_2,e_3+ie_6)$ respectively.  If we
then apply \eqref{3c9} and \eqref{3c10} to these vectors at their
corresponding points, we obtain
\bea
ge_1\times ge_2 =ge_4\times ge_7 =ge_6\times ge_5 \label{3c11} \\
ge_1\times ge_3 =ge_6\times ge_4 =ge_7\times ge_5 \label{3c12} \\
ge_1\times ge_6 =ge_5\times ge_2 =ge_4\times ge_3 \label{3c13} \\
ge_1\times ge_7 =ge_5\times ge_3 =ge_2\times ge_4 \label{3c14}
\eea
The vectors $\{ ge_1,..,ge_7\}$ are orthonormal with respect
to the Euclidean product $(,)$ in $\cc ^7$ since
$$
(ge_i,ge_j)=\ip{ge_i}{\ol{g}e_j}=\ip{g^tge_i}{e_j}=\dt_{ij}.
$$
Recalling that $(*\times *,*)$ is skew-symetric, we see that
$$
(ge_i\times ge_j,ge_j)=0.
$$
Thus it follows from \eqref{3c11} that $ge_1\times ge_2 =\pm
ge_3$.

If $ge_1\times ge_2 =ge_3$ then we can use directly \eqref{3c12},
\eqref{3c13} and \eqref{3c14} to show that $\{ ge_1,..ge_7\}$ is a
$G_2$-basis for $\cc ^7$ and hence $g\in G_2^{\cc}$.

Similarly, if $ge_1\times ge_2 =-ge_3$ then we can repeat the same
process above to deduce that $-g \in G_2^{\cc}$. \epf

\section{Superminimal Surfaces in $S^6$}

\begin{deff} \label{stsu2}
Let $S$ be a Riemann surface.  We say that a smooth map $f:S\to
S^6$ is an {\bf almost complex curve} of the nearly K\" ahler
$S^6$ if $f_*$ is complex linear.
\end{deff}
Therefore, using a local complex coordinate $z = x+iy$ for $S$ we
can characterize these curves by

\bea
\partial _yf=f\times \partial _xf. \label{su8}
\eea

It follows that almost complex curves of $S^6$ are weakly
conformal maps which are also harmonic maps because if we differentiate again
\eqref{su8}, we obtain $f\times (\partial_{xx}f+\partial
_{yy}f)=0$, yielding that $\partial_{xx}f + \partial_{yy}f$ is normal to $S^6$,
and hence $f$ is harmonic in accordance with \cite{EL}.

Therefore an almost complex curve $f:S \to S^6$ determines a
harmonic sequence of maps $\psi _k:S\to \cc P^6$ so that $\psi
_0=[f]$.  Using this sequence and some invariants associated to
their elements, a full classification of these curves was obtained
in \cite{bvw1} according to the following four types:
\begin{en3}
\item Linearly full in $S^6$ and superminimal, \item Linearly full
in $S^6$ but not superminimal, \item Linearly full in a totally
geodesic $S^5$ in $S^6$, \item Totally geodesic.
\end{en3}
A result of Bryant \cite{bryant1} highlights the importance of the
Type-I almost complex curves of $S^6$.  Bryant has shown that
every compact Riemann surface of any genus can be realised as such
an almost complex curve of the 6-sphere.

In this section we shall obtain explicitly all the 0-point and
2-point ramified linearly full almost complex 2-spheres of the
6-sphere. This is done by using the normal form for such surfaces
as given by Theorem \eqref{stsu5}.

In particular, we will also prove that these surfaces are uniquely
determined by their singularity type up to $G_2^{\cc
}$-equivalence of their directrix curves. It is worthwhile
mentioning here that a similar result in the more general
situation of harmonic 2-spheres of $S^n$ and $\cc P^n$ has been
obtained in \cite{bw3} but replacing, of course, the group
$G_2^{\cc}$ by $SO(n+1,\cc )$.

\bpro \label{stsu6}

Let $f:S\to S^6$ be a
linearly full superminimal almost complex curve. If $\big( \psi
_j=[f_j] \big) _{j=0}^{6}$ is the harmonic sequence corresponding
to the harmonic map $\psi _3=[f_3]=[f]$ then the meromorphic local
sections $f_k:\cc \to \cc ^7$ have the following multiplication
table for $f_i \times f_j$, where the cross product $\times$ is
extended $\cc$-linearly to $\cc^7$:
\begin{eqnarray}
\begin{array}{|c||c|c|c|c|c|c|c|}
\hline
i\backslash j & \hs 0 &\hs 1 &\hs 2 &\hs 3 &\hs 4 &\hs 5 &\hs 6 \\
\hline \hline 0&0&0&0&-if_0&-2if_1&-2if_2&-if_3   \\ \hline
1&0&0&if_0&if_1&0&-if_3&-if_4       \\ \hline
2&0&-if_0&0&if_2&if_3&0&-if_5       \\ \hline
3&if_0&-if_1&-if_2&0&if_4&if_5&-if_6    \\ \hline
4&2if_1&0&-if_3&-if_4&0&2if_6&0         \\ \hline
5&2if_2&if_3&0&-if_5&-2if_6&0&0         \\ \hline
6&if_3&if_4&if_5&if_6&0&0&0 \\ \hline
\end{array}
\label{table-f}
\end{eqnarray}
Furthemore, the following relation holds \bea
|f_4|^2|f_5|^2=2|f_6|^2. \label{su10} \eea

\epro

The linearly fullness and superminimality conditions can easily be used to
calculate as in \cite{bvw1} both the condition \eqref{su10} and the following products:
$f_3\times f_4=if_4$, $f_3\times f_5=if_5$ and $f_3\times f_6=-if_6$.
Then all the other products can be obtained by taking derivatives $\fr{\partial}{\partial z}$
and $\fr{\partial}{\partial \ol z}$ and the properties given by formulae
\eqref{su3}, \eqref{su4} and \eqref{su5a}.

\brem

It is worth mentioning that the condition \eqref{su10}
characterizes the linearly full superminimal almost complex curves
of the 6-sphere (cf. \cite{bvw1}) in the sense that a weakly
conformal harmonic map $f:S \to S^6$ is $O(7)$-congruent to a
linearly full superminimal almost complex curve if and only if
\eqref{su10} holds.

\erem

We say that a map $\psi$ from a Riemann
surface $S$ into $Q^5$ is superhorizontal if at each point of $S$,
$\psi _*$ takes values in the superhorizontal distribution.  We
shall recall now the 1-1 correspondence (cf.\cite{bw4})
between superminimal almost complex curves in $S^6$ and
holomorphic superhorizontal curves in $Q^5$. We intend to make use
of this correspondence later on in this article in order to work
out explicit examples of superminimal 2-spheres of $S^6$.

By using \eqref{3c5}, the superhorizontal condition of $\psi =[f]$
can be described analitically as follows.
\begin{eqnarray*}
\mbox{$\psi$ is superhorizontal} & \Llra &f_*|_p(T_pS) \subseteq
H'_{\psi (p)}
= \{v \in T_{[f(p)]}Q^5 :f\times v = 0\} \\
& \Llra &f \times f_*(a\frac{\partial}{\partial z} +
b\frac{\partial}{\partial \bar{z}})
= 0 \hs \hs \forall \hs a,b \in \rr \\
& \Llra & f\times f_z = 0 \hs \hs \mbox{ and } \hs \hs f \times
f_{\bar{z}} = 0.
\end{eqnarray*}
thus a holomorphic map $\psi :S\to Q^5$ is superhorizontal if and
only if \bea f\times \fr{df}{dz} = 0.\label{su20} \eea As a
consequence of this characterization, The theorem below is proved by Bryant in \cite{bryant1} and Bolton at al have
given a tidy treatment in \cite{bw4}.

\begin{thm} \label{stsu7}
A map $g: S \to Q^5$ is linearly full, holomorphic and superhorizontal if and
only if $\psi =\pi (g): S \to S^6$ is a linearly full superminimal almost
complex curve in $S^6$ with directrix curve $g$ where $\pi$
denotes the twistor map from $Q^5$ onto $S^6$.
\end{thm}

Here we can prove the following characterization for almost
complex curves in $S^6$ in terms of their directrix curves.

\bthm \label{stsu8} Let $f$ and $\td f$ be linearly full almost
complex 2-spheres of $S^6$. Then their directrix curves are
projectively equivalent if and only if they are also
$G_2^{\cc}$-equivalent. \ethm \bpf \come The converse in the
theorem is obvious since $G_2^\cc$ is a subgroup of $SO(7,\cc )$.

\go Let $(\psi _j=[f_j])_{j=0}^6$ denote the harmonic sequence
corresponding to the harmonic map $[f]$ and let $\wt{\psi}_0$
denote the directrix curve of the harmonic map $[\td f]$. By
assumption there exists an element $[A]\in PGL(7,\cc )$ such that
$\wt{\psi}_0=[Af_0]$.

It is shown in \cite{bw3} (Theorem 3.3) that two linearly full
harmonic 2-spheres of $S^{2m}$ are projectively equivalent if and
only if they are $SO(2m+1, \cc )$-equivalent. Thus we can assume
in our particular case here that $A$ lies in $SO(7, \cc )$.

According to theorem \eqref{stsu7} the map $Af_0:S^2 \to \cc ^7$
is holomorphic and superhorizontal and hence \bea Af_0\times
Af'_0=0. \mbox{\hs \hs $(f'_0=\fr{df_0}{dz})$} \label{su21} \eea
We shall make use in the sequel of the following properties
satisfied by the functions $\big( \al _j \big)_{j=0}^6$ defined by
equation \eqref{su3}. \bea \al _3=&\zzz \zzz 0.
&\mbox{\hs Follows from \eqref{su5b},} \notag \\
\al _j=&\zzz -\al _{6-j}.
&\mbox{\hs Follows from \eqref{su5a},} \\
\al _6=&\al _5 + \al _4. &\mbox{\hs Follows from \eqref{su10}.}
\eea We differentiate \eqref{su21} with respect to $z$ and use
\eqref{su3}, obtaining in this way the cross product between
different vectors $Af_j$. By repeating this process we can derive
some relations among the cross product of the vectors $Af_j$,
namely
\bea
Af_0\times Af_1 &=&0
\label{su22} \\
Af_0\times Af_2 &=&0
\label{su23} \\
Af_0\times Af_3 &=&-Af_1\times Af_2
\label{su24} \\
Af_0\times Af_4 &=&-2Af_1\times Af_3
\label{su25} \\
Af_0\times Af_5 &=&-2Af_2\times Af_3 -3Af_1\times Af_4
\label{su26} \\
Af_0\times Af_6 &=& -5Af_2\times Af_4 -4Af_1\times Af_5
+3\al_5Af_1\times Af_4. \label{su27} \eea As $A\in SO(7,\cc )$, it
follows from \eqref{su5} that \bea (Af_i,Af_j)
&=&\ip{Af_i}{\ol{A}\ol{f}_j}
\notag \\
&=&\ip{A^tAf_i}{\ol{f}_j}
\notag \\
&=&\ip{f_i}{\ol{f}_j}
\notag \\
&=&(-1)^i\dt_{i,6-j}. \label{su28} \eea The vectors $\{
Af_0,\hdots ,Af_6\}$ form a basis for $\cc ^7$ since $A\in
SO(7,\cc )$. Thus, from \eqref{su24} we see immediately that
$Af_0\times Af_3$ can be written as the linear combination:
$$
Af_0\times Af_3 =aAf_0+bAf_1+cAf_2.
$$
But if we take the cross product of this equation with $Af_2$ and
$Af_1$, then it follows that $b=c=0$. Indeed,
\begin{eqnarray*}
\begin{array}{rlll}
0&=&Af_2\times (Af_0\times Af_3) &\mbox{use \eqref{su28}}
\\
&=&bAf_2\times Af_1 &\mbox{use \eqref{su23}}
\\
&=&b(Af_0\times Af_3), &\mbox{use \eqref{su24}}
\\
&\mbox{and}&
\\
0&=&Af_1\times (Af_0\times Af_3) &\mbox{use \eqref{su28}}
\\
&=&cAf_1\times Af_2 &\mbox{use \eqref{su22}}
\\
&=&-c(Af_0\times Af_3). &\mbox{use \eqref{su24}}
\end{array}
\end{eqnarray*}
On the other hand, if we also take the cross product with $Af_6$,
we see that $a=\pm i$. Indeed,
\begin{eqnarray*}
\begin{array}{rll}
a(Af_0\times Af_6)
&=&-Af_6\times (Af_0\times Af_3) \\
&=&(Af_6\times Af_0)\times Af_3 - Af_3,
\end{array}
\end{eqnarray*}
and the Euclidean product of this with $Af_0$ gives:
\begin{eqnarray*}
\begin{array}{rll}
-1&=&a(Af_0\times Af_6,Af_3)\\
&=&-a(Af_0\times Af_3,Af_6) \\
&=&-a^2(Af_0,Af_6) \\
&=&a^2.
\end{array}
\end{eqnarray*}
Let us first assume the case $Af_0\times Af_3=-iAf_0$. Then
\eqref{su25} yields $Af_1\times Af_3=iAf_1$, indeed
$$
\begin{array}{lll}
-2Af_1\times Af_3&=&Af_0\times Af_4 \\
&=&i(Af_0\times Af_3)\times Af_4 \\
&=&-i(Af_0\times Af_4)\times Af_3 \\
&=&2i(Af_1\times Af_3)\times Af_3 \\
&=&-2iAf_1.
\end{array}
$$
Thus from \eqref{su25} we get
$$
Af_0\times Af_4=-2iAf_1.
$$
And this yields:
$$
Af_1\times Af_4=\fr{i}{2}(Af_0\times Af_4)\times Af_4=0. \hs \hs
\mbox{(Using \eqref{su8}).}
$$
Now, by using the equations \eqref{su21},...,\eqref{su28}, we can
carry on with this process to determine all the cross products of
the vectors $\{ Af_0,..,Af_6\}$ and to verify that they satisfy
the multiplication table \eqref{table-f} in the following sense
$$
A(f_i\times f_j)=Af_i\times Af_j.
$$
Therefore, $A$ is an element of $G_2^{\cc }$ since $\{
f_0,..,f_6\}$ is a basis for $\cc ^7$.

In the other case to be considered, that is, when $Af_0\times
Af_3=iAf_0$, we can use the same procedure as above to prove that
$-A\in G_2^{\cc }$. However, this contradicts our assumption that
$A\in SO(7,\cc )$.\epf

\vspace{3mm}

Let $f:S^2\to S^6$ be a $k$-point ramified $(k\leq 2)$ linearly full
almost complex curve and let $\psi _j=[f_j]$ $(j=0,\hdots ,6)$
denote the harmonic sequence corresponding to the harmonic map
$\psi _3=[f_3]=[f]$. Then according to Theorem \eqref{stsu5}, we can
find a local complex coordinate $z$ for $S^2$ and a basis $\{
v_0,..,v_6 \}$ for $\cc ^7$, so that the directrix curve $\psi
_0=[f_0]$ can be expressed by:
\begin{align}
\begin{split}
f_0=
&v_0+z^{k_1}v_1+z^{k_1+k_2}v_2+z^{k_1+k_2+k_3}v_3+ \\
&+z^{k_1+k_2+2k_3}v_4+z^{k_1+2k_2+2k_3}v_5+z^{2k_1+2k_2+2k_3}v_6.
\label{su29}
\end{split}
\end{align}

Using the meromorphic sections $f_j$ we can choose a particular
orthonormal basis $\{ u_0,\hdots ,u_6\}$ for $\cc ^7$ so that each
vector $u_j$ spans the same complex line bundle as $f_j$. Indeed,
as for $z\not= 0$ the function $\fr{f_j}{|f_j|}$ takes values in
the sphere $S^{13}$ then for each $j\in\{0,\hdots ,6\}$, by compactness,
there exists a sequence $w_{jk}\mapsto 0$ so that  we have
\begin{equation}
\lim_{k\mapsto \infty}{\fr{f_j}{|f_j|}(w_{jk})}=u_j \label{su30}
\end{equation}
Now, we can notice that the vectors $u_j$ are not obtained by orthonormalization
but the triangular display below appears naturally from manipulation of equations \eqref{su29}
and \eqref{su30}. For example, the second equation is obtained in the following way:
\begin{equation} \label{su35}
  f_1+\alpha _0f_0 = \fr{\partial f_0}{\partial z}=k_1z^{k_1-1}v_1+\hdots
\end{equation}
Consequently, we obtain a vector function $\rho (z)$ satisfying $\lim _{z\to 0}\rho (z)=0$ and
\begin{equation*}
  |f_1(w_{jk})|u_1+\alpha _0|f_0(w_{jk})|u_0 =k_1(w_{jk})^{k_1-1}.[v_1+\rho (w_{jk})]
\end{equation*}
and hence $v_1$ is a linear combination of the vectors $u_1$ and $u_0$. Similarly,
by differentiating successively equation \eqref{su35}, we obtain the following
triangular shape:

\vspace{2mm}
\noindent
$ v_0 =a_{\sss (0,0)}u_0 $
\\
$ v_1 =a_{\sss (1,0)}u_0+a_{\sss (1,1)}u_1 $
\\
$ v_2 =a_{\sss (2,0)}u_0+a_{\sss (2,1)}u_1+a_{\sss (2,2)}u_2 $
\\
$ v_3 =a_{\sss (3,0)}u_0+a_{\sss (3,1)}u_1+a_{\sss (3,2)}u_2
+a_{\sss (3,3)}u_3 $
\\
$ v_4 =a_{\sss (4,0)}u_0+a_{\sss (4,1)}u_1+a_{\sss (4,2)}u_2
+a_{\sss (4,3)}u_3+a_{\sss (4,4)}u_4 $
\\
$ v_5 =a_{\sss (5,0)}u_0+a_{\sss (5,1)}u_1+a_{\sss (5,2)}u_2
+a_{\sss (5,3)}u_3+a_{\sss (5,4)}u_4+a_{\sss (5,5)}u_5 $
\\
$ v_6 =a_{\sss (6,0)}u_0+a_{\sss (6,1)}u_1+a_{\sss (6,2)}u_2
+a_{\sss (6,3)}u_3+a_{\sss (6,4)}u_4+a_{\sss (6,5)}u_5+a_{\sss
(6,6)}u_6$

\vspace{2mm}

\noindent where the scalars $(a_{\sss (i,j)})$ appearing in the linear
combinations are complex numbers.

From \eqref{su30} and Proposition \eqref{stsu6} we can easily
determine the cross product of the vectors $u_j$ and consequently
also of the vectors $v_j$. Namely, the vectors $u_j$ have the
following multiplication table for $u_i \times u_j$:
\begin{eqnarray*}
\wc
\begin{array}{|c||c|c|c|c|c|c|c|}
\hline
i\backslash j & \hs 0 &\hs 1 &\hs 2 &\hs 3 &\hs 4 &\hs 5 &\hs 6 \\
\hline \hline 0&0&0&0&-iu_0&-i\sqrt{2}u_1&-i\sqrt{2}u_2&-iu_3 \\
\hline 1&0&0&i\sqrt{2}u_0&iu_1&0&-iu_3&-i\sqrt{2}u_4   \\ \hline
2&0&-i\sqrt{2}u_0&0&iu_2&iu_3&0&-i\sqrt{2}u_5   \\ \hline
3&iu_0&-iu_1&-iu_2&0&iu_4&iu_5&-iu_6    \\ \hline
4&i\sqrt{2}u_1&0&-iu_3&-iu_4&0&i\sqrt{2}u_6&0       \\ \hline
5&i\sqrt{2}u_2&iu_3&0&-iu_5&-i\sqrt{2}u_6&0&0       \\ \hline
6&iu_3&i\sqrt{2}u_4&i\sqrt{2}u_5&iu_6&0&0&0 \\ \hline
\end{array}
\label{table-u}
\end{eqnarray*}

Now, we notice that the coefficient $a_{\sss (6,6)}$ must be
non-zero since $f$ is linearly full and also $f_0$ is a
holomorphic superhorizontal curve because of the characterization
given in Theorem \eqref{stsu7}. These facts together with equation
\eqref{su20}, give us a cumbersome but straigthforward calculation
to determine the following example.

\bins \label{example-1} Let $\{ e_1,\hdots ,e_7\}$ denote a
$G_2$-basis for $\rr ^7$. Let $k_1,k_2$ denote positive integers.
Then the holomorphic map $\psi _0 =[f_0]:S^2\to \cc P^6$,
determined by the polynomial $f_0(z)=\sum_{j=1}^{7}a_j(z)e_j$
where the $a_j(z)$ are given by

\vspace{2mm}
\noindent
$ a_1(z)=
-\fr{3\sqrt{30}k_2(k_1+k_2)}{(3k_1+k_2)(2k_1+k_2)}z^{k_1+k_2}
+\fr{\sqrt{30}}{2}z^{3k_1+k_2}, $
\\
$ a_2(z)= \fr{15\sqrt{3}k_1k_2}{(3k_1+2k_2)(2k_1+k_2)}z^{k_1}
+\sqrt{3}z^{3k_1+2k_2}, $
\\
$ a_3(z)=
\fr{i45\sqrt{2}k_1k_2^2(k_1+k_2)}{(3k_1+2k_2)(3k_1+k_2)(2k_1+k_2)^2}
+i\fr{\sqrt{2}}{2}z^{4k_1+2k_2}, $
\\
$ a_4(z)= \fr{6\sqrt{5}k_2}{2k_1+k_2}z^{2k_1+k_2}, $
\\
$ a_5(z)=
i\fr{3\sqrt{30}k_2(k_1+k_2)}{(3k_1+k_2)(2k_1+k_2)}z^{k_1+k_2}
+i\fr{\sqrt{30}}{2}z^{3k_1+k_2}, $
\\
$ a_6(z)= -i\fr{15\sqrt{3}k_1k_2}{(3k_1+2k_2)(2k_1+k_2)}z^{k_1}
+i\sqrt{3}z^{3k_1+2k_2}, $
\\
$ a_7(z)=
-\fr{45\sqrt{2}k_1k_2^2(k_1+k_2)}{(3k_1+2k_2)(3k_1+k_2)(2k_1+k_2)^2}
+\fr{\sqrt{2}}{2}z^{4k_1+2k_2}, $

\vspace{2mm}
\noindent is the directrix curve of a linearly full $S^1$-symmetric almost
complex $2$-sphere in $S^6$ with the same singularity type
$(k_1-1,k_2-1,k_1-1,k_1-1,k_2-1,k_1-1)$ at $z=0$ and $z=\infty$.
\eins

\begin{thm} \label{k-ram}
Let $f:S^2\to S^6$ be a k-point ramified $(k\leq 2)$ linearly full
almost complex curve with any singularities at $z=0$ and $z=\infty$.
Then for a suitable choice of complex coordinate on $S^2$, the harmonic map
$[f(z)]:S^2 \to \cc P^6$ has the same singularity type
$(k_1-1,k_2-1,k_1-1,k_1-1,k_2-1,k_1-1)$ at $z=0$ and $z=\infty$.
Moreover, the directrix curve of $f$ is $G_2^{\cc}$-equivalent to
the $S^1$-symmetric curve given in the Example \eqref{example-1}.
\end{thm}
\bpf
Let $\psi _0:S^2\to Q^5$ given by $\psi_0(z)=[f_0(z)]$
be the directrix curve of the map
$[f_3(z)]=[f(z)]$. The first part of the statement follows from item (ii)
of Theorem \eqref{stsu5} and the following observation.

By comparing the exponents of the variable $z$ appearing in the
polynomial $f\times f_z=0$, we obtain the symmetry $k_3=k_1$ in
the singularity type.

Now, we can use the condition given by the equation \eqref{su20}
together with the multiplication table for the products $u_i\times u_j$
in order to determine the vectors $(v_j)$ in a more simplified
way in terms of the vectors $u_j$, and expressed only in terms
of the 8 complex parameters $r_1= a_{(4,4)}, \hs r_2= a_{(5,4)}, \hs
r_3 = a_{(6,3)}, \hs r_4 = a_{(6,4)},\hs r_5= a_{(6,5)}, \hs r_6= a_{(6,2)},
\hs r_7 = a_{(6,1)}, \hs r_8 = a_{(5,5)}$, in the following way:

\noindent
$ v_0=
{\frac{k_1k_2^2(k_1+k_2)r_1^2r_8^2}{(3k_1+2k_2)(3k_1+k_2)(2k_1+k_2)^2}}
u_0, $
\\
$ v_1= \frac{k_1k_2r_1^2r_8}{(3k_1+2k_2)(2k_1+k_2)} (r_5u_0+u_1),
$
\\
$ v_2= \frac{k_2(k_1+k_2)r_1r_8}{(2k_1+k_2)(3k_1+k_2)}
((r_2r_5-r_4r_8)u_0+r_2u_1+r_8u_2), $
\\
$ v_3= \frac{k_2r_1r_8}{(2k_1+k_2)}
(\sqrt{2}r_3u_0+2r_4u_1+2r_5u_2+\sqrt{2}u_3), $
\\
$ v_4= \frac{r_1}{2}
((\sqrt{2}r_3r_5-2r_6)u_0+(2r_4r_5-\sqrt{2}r_3)u_1
+2r_5^2u_2+2\sqrt{2}r_5u_3+2u_4), $
\\
$ v_5=
(\fr{\sqrt{2}}{2}r_2r_3r_5-\fr{\sqrt{2}}{2}r_3r_4r_8+r_7r_8-r_2r_6)u_0
+(r_2r_4r_5-\fr{\sqrt{2}}{2}r_2r_3-r_4^2r_8)u_1, $
\\
$ \phantom{v_2=} +(r_2r_5^2-\fr{\sqrt{2}}{2}r_3r_8-r_4r_5r_8)u_2
+\sqrt{2}(r_2r_5-r_4r_8)u_3+r_2u_4+r_8u_5, $
\\
$ v_6= (r_5r_7+\fr{1}{2}r_3^2-r_4r_6)u_0
+r_7u_1+r_6u_2+r_3u_3+r_4u_4+r_5u_5+u_6. $

\vspace{2mm}

By Corollary \eqref{stsu5b} we can assume $f$ to be
$S^1$-symmetric. Theorem \eqref{stsu5} shows that the
$S^1$-symmetric linearly full almost complex 2-spheres are
characterized by the orthogonality of the vectors $(v_j)$ and
hence according to the formulae above we must have $r_2=\hdots
=r_7=0$. Thus, the directrix curve is described by the 2-parameter
family
\begin{equation*}
\begin{split}
f_0(z)=
&{\frac{k_1k_2^2(k_1+k_2)r_1^2r_8^2}{(3k_1+2k_2)(3k_1+k_2)(2k_1+k_2)^2}}
u_0 +\frac{k_1k_2r_1^2r_8}{(3k_1+2k_2)(2k_1+k_2)} z^{k_1}u_1
\\
&+\frac{k_2(k_1+k_2)r_1r_8^2}{(2k_1+k_2)(3k_1+k_2)} z^{k_1+k_2}u_2
+\sqrt{2}\frac{k_2r_1r_8}{(2k_1+k_2)} z^{2k_1+k_2}u_3
\\
&+r_1z^{3k_1+k_2}u_4 +r_8z^{3k_1+2k_2}u_5 +z^{4k_1+2k_2}u_6.
\end{split}
\end{equation*}
Now, we shall apply a suitable conformal transformation to the
domain and also apply an appropriate element of $G_2$ to the
co-domain in order to prove that $f$ is indeed equivalent to the
curve given in the example above.

Let $r$ be a complex root for the equation \bea
r^{2k_1+k_2}r_1r_8=\sqrt{90}. \label{su31} \eea Then we shall
consider the conformal transformation $z \mapsto rz$, and the
element $A\in G_2$ defined by
$$
\begin{array}{lll}
Au_0:=(\fr{90}{r_1^2r_8^2})u_0, \hs
&Au_1:=(\fr{15\sqrt{6}}{r^{k_1}r_1^2r_8})u_1, \hs
&Au_2:=(\fr{6\sqrt{15}}{r^{(k_1+k_2)}r_1r_8^2})u_2, \hs
\\
Au_3:=(\fr{\sqrt{90}}{r^{(2k_1+k_2)}r_1r_8})u_3, \hs
&Au_4:=(\fr{\sqrt{15}}{r^{(3k_1+k_2)}r_1})u_4,  \hs
&Au_5:=(\fr{\sqrt{6}}{r^{(3k_1+2k_2)}r_8})u_5, \hs
\\
Au_6:=(\fr{1}{r^{(4k_1+2k_2)}})u_6. &&
\end{array}
$$
Using \eqref{su31} and the multiplication table \eqref{table-u}
for the vectors $u_j$ we deduce that
$$
A(u_i\times u_j)=Au_i\times Au_j,
$$
which implies that $A\in G_2$. Thus, the holomorphic curve $Af_0$
is reduced to
\begin{equation*}
\begin{split}
Af_0(z)=
&{\frac{90k_1k_2^2(k_1+k_2)}{(3k_1+2k_2)(3k_1+k_2)(2k_1+k_2)^2}}
u_0 +\frac{15\sqrt{6}k_1k_2}{(3k_1+2k_2)(2k_1+k_2)} z^{k_1}u_1
\\
&+\frac{6\sqrt{15}k_2(k_1+k_2)}{(2k_1+k_2)(3k_1+k_2)}
z^{k_1+k_2}u_2
+\sqrt{2}\frac{6\sqrt{5}k_2}{(2k_1+k_2)}z^{2k_1+k_2}u_3
\\
&+\sqrt{15}z^{3k_1+k_2}u_4 +\sqrt{6}z^{3k_1+2k_2}u_5
+z^{4k_1+2k_2}u_6.
\end{split}
\end{equation*}
Using again that multiplication table we can also deduce by
straigthforward calculations that the vectors $e_j\in \rr ^7
\hs(j=1,\hdots ,7)$ defined by
\\
$ u_0=\fr{1}{\sqrt{2}}(-e_7 +ie_3 ), $
\hspace{2mm}
$ u_1=\fr{1}{\sqrt{2}}(e_2 -ie_6 ), $
\hspace{2mm}
$ u_2=\fr{1}{\sqrt{2}}(-e_1 +ie_5 ), $
\hspace{2mm}
$ u_3=e_4, $
\\
$ u_4=\fr{1}{\sqrt{2}}(e_1 +ie_5 ), $
\hspace{2mm}
$ u_5=\fr{1}{\sqrt{2}}(e_2 +ie_6 ), $
\hspace{2mm}
$ u_6=\fr{1}{\sqrt{2}}(e_7 +ie_3 ), $
\\
form a $G_2$-basis for $\rr ^7$ and the holomorphic curve $Af_0(z)$
is written in terms of this basis exactly as the one we gave in
the example.\epf

\bcor \label{moduli-1}
If $H^{\sss 0,0}$ is the space of linearly full totally unramified
almost complex maps of $S^2$ into $S^6$ then $\C H^{\sss 0,0}=G_2^{\cc}$.
\ecor
\bpf
Indeed, this follows from the theorem above and the
fact that the harmonic sequence corresponding to a harmonic map
$[f]$, where $f\in \C H^{\sss 0,0}$, is uniquely determined by its
directrix curve. Some care is required here since the composition
of $f$ with the antipodal map of $S^6$ gives also a harmonic map
with that same directrix curve. However, the map $-f$ is
fortunately an almost anticomplex curve as we can see from
\eqref{su8}.\epf

\vspace{2mm}

Let $M$ denote the quotient set of the manifold $N=\{ (p,q)\in
S^2\times S^2 / p\not= q\}$ by the equivalence relation:
$(p,q)\cong (a,b)$ if and only if $p=b$ and $q=a$.

\bcor
\label{moduli-2} Let $H^{\sss r_1,r_2}$ denote the space of
linearly full almost complex maps of $S^2$ into $S^6$ with 2
higher singularities each of type $(r_1,r_2,r_1,r_1,r_2,r_1)$.
Then $\C H^{\sss r_1,r_2}=M\times G_2^{\cc}$.
\ecor

Let $ \psi _0,\hdots , \psi _n : S^2 \to \cc P^n $ be a harmonic
sequence with corresponding local lifts $f_0,\hdots ,f_n:
S^2\backslash W \to \cc P^n$ given in accordance with \eqref{su3}
and \eqref{su4}, where $W$ is the set of all singularities of the
harmonic maps $\psi _p$. Bolton et al have proved in \cite{bjrw}
that when $\psi _p$ is an immersion, the area $A(\psi _p)$ of
$S^2$ with the metric induced by $\psi _p$ is given by \bea A(\psi
_p)=\pi (\dt _{p-1}+\dt _p), \label{area} \eea where $\dt _{-1}=0$
and $\dt _p$ is the degree of the $(p-1)$-st osculating curve $\sg
_{p-1}$. Moreover, they calculate this degree in terms of the $\gm
_p$ invariants. Namely, \bea \dt _p=\fr{1}{2\pi i}\int _{S^2}\gm
_pd\ol z\wedge dz. \label{delta-1} \eea Bolton et al carry on
working out the following global Pl\" ucker formula, relating the
ramification indices $R_p$ of the curves $\sg _{p-1}$ to the
degrees $\dt _p$ by \bea R_p=-2-\dt _{p-2}+2\dt _{p-1}-\dt _{p}.
\label{gpf} \eea Finally, they also write down the degrees $\dt
_p$ in terms of the the ramification indices $R_p$ as follows \bea
\dt _p= (p+1)(n-p)+ \fr{n-p}{n+1} \sum _{k=0}^{p-1}(k+1)R_k+
\fr{p+1}{n+1} \sum _{k=p}^{n-1}(n-k)R_k. \label{delta-2} \eea
Using these results for the case $n=6$ we can now produce the
following consequence \blem Let $f:S^2 \to S^6$ be a linearly full
almost complex curve. Let $(\psi _p)_{p=0}^{6}$ be the harmonic
sequence determined by $f$. Then the ramification indices $R_p$ of
the associated osculating curves $\sg _{p-1}$ satisfy
\begin{en1}
\item $R_j=R_{7-j}$\hs for \hs $j=1,\hdots ,6$. \item $R_3=R_2$.
\end{en1}
\elem \bpf Considering that $f$ is superminimal (see paragraph
after Definition \eqref{stsu1}) item (i) follows from direct
application of \eqref{gpf}, \eqref{delta-1},\eqref{su5a} and
\eqref{su10}, while item (ii) follows from \eqref{gpf},
\eqref{delta-1} and \eqref{su10}. \epf
\bpro
Let $f:S^2 \to S^6$ be a linearly full almost complex curve with total singularity
type $(R_1,\hdots ,R_6)$. Then the area $A(\psi )$ of the harmonic
map $\psi =[f]:S^2\to \cc P^6$ is given by \bea A(\psi )=4\pi
(6+2R_1 +R_2). \label{area-2} \eea \epro \bpf Using the lemma
above and \eqref{delta-2} we have
$$\dt _3=12+4R_1+2R_2=\dt _4.$$
Thus, the Corollary follows from \eqref{area}.\epf

\bthm Let $\C H ^d$ be the space of linearly full almost complex maps of $S^2$
into $S^6$ of area $4\pi d$. Then $d\ge 6$ and
\begin{en1}
\item $\C H^6=\C H^{\sss 0,0}=G_2^{\cc}$, \item $\C H^7$ is empty,
\item $\C H^8=\C H^{\sss (0,1)}=M\times G_2^{\cc}$.
\end{en1}
Furthermore, every element of $\C H^8$ has directrix curve
$G_2^{\cc}$-equivalent to the following $S^1$-symmetric case \bea
f(z) && = (70\sqrt{15}z^5-126\sqrt{15}z^3)e_1
+(70\sqrt{6}z^7+75\sqrt{6}z)e_2
\notag \\
&& +(i135+i70z^8)e_3 +210\sqrt{10}z^4e_4
+(i70\sqrt{15}z^5+i126\sqrt{15}z^3)e_5
\notag \\
&& +(i70\sqrt{6}z^7-i75\sqrt{6}z)e_6 +(-135+70z^8)e_7.
\label{lowest} \eea
\ethm
\bpf Item (i) and $d\ge 6$ follow are
consequence of \eqref{area-2}, whilst item (ii) follows from the
fact already noted in the Remark \eqref{stsu5a} that there does
not exist a 1-point ramified harmonic 2-sphere in $\cc P^n$.

Now to prove (iii) we start noting that according to
\eqref{area-2} the first possibility for the singularity type is
$R_0=0$ and $R_1=2$. This case implies that $[f]$ is 2-point
ramified with singularity type at each point given by
$(0,1,0,0,1,0)$ and hence (iii) and \eqref{lowest} follow from
Theorem \eqref{k-ram}. \epf

\addcontentsline{toc}{chapter}{\numberline{}References.}


\begin{thebibliography}{99}
\bibitem{barbosa1}
Barbosa J.L.M.:
On minimal immersions of $S^2$ into $S^{2m}$.
Trans. Amer. Math. Soc. 210 (1975), 75-106.
\bibitem{bw4}
Bolton J., Woodward L.M.
Special submanifolds of $S^6$   with its $G_2$ geometry.
Geometry, topology and physics (Campinas,1996), 59-68, de Gruyter, Berlin, 1997.
\bibitem{bw1}
Bolton J., Woodward L.M.:
Congruence theorems for harmonic maps from a Riemann surface
into $\cc P^n$ and $S^n$.
J. London Math. Soc. (2) 45 (1992), 363-376.
\bibitem{bw2}
Bolton J., Woodward L.M.:
On the Simon conjecture for minimal immersions with $S^1$-symmetry.
Math. Z. 200 (1988), 111-121.
\bibitem{bw3}
Bolton J., Woodward L.M.:
Moduli spaces of harmonic 2-spheres.
Geometry and Topology of Submanifolds IV, World Scientific (1992)
143-151
\bibitem{bvw1}
Bolton J., L.Vrancken, L.M.Woodward:
On Almost Complex Curves in the Nearly K\" ahler 6-Sphere.
Quart. J. Math. Oxford (2) 45 (1994), 407-427.
\bibitem{bvw2}
Bolton J., L.Vrancken, L.M.Woodward:
Minimal immersions of $S^2$ and $\rr P^2$ into $\cc P^n$
with few higher order singularities.
Math. Proc. Camb. Phil. Soc. 111 (1992), 93-101.
\bibitem{bjrw}
Bolton J., Jensen G.R., Rigoli M., Woodward L.M.:
On conformal minimal immersions of $S^2$ into $\cc P^n$.
Math. Ann. 279 (1988) 599-620.
\bibitem{bryant1}
Bryant R.L.:
Submanifolds and special structures on the octonians.
J. Diff. Geom. 17 (1982) 185-232.
\bibitem{dz1}
Din A.M., Zakrzewski W.J.:
General classical solutions in $\cc P^{n-1}$-model.
Nuclear Phys. B 174 (1980), 397-406.
\bibitem{EL}
Eells, J.; Lemaire, L.: 
A report on harmonic maps. 
Bull. London Math. Soc. 10 (1978), no.1, 1-68.
\bibitem{ew1}
Eells J., Wood J.C.:
Harmonic maps from surfaces to complex projective spaces.
Adv. Math. 49 (1983) 217-263.
\bibitem{gh}
Griffiths P., Harris J.:
Principles of algebraic geometry.
London-New York, Wiley 1978.
\bibitem{martins1}
Martins J.K.: Congruence of Hypersurfaces in $S^6$ and in $CP^n$.
Bol. Soc. Bras. Mat. Vol 32 (2001), n.1, 83-105.
\bibitem{martins2}
Martins J.K.: Hopf Hypersurfaces in Space Forms. Bull. Braz. Math.
Soc. New Series Vol. 35 (2004), n.3, 453-472.
\bibitem{wolfson1}
Wolfson J.G.:
Harmonic sequences and harmonic maps of surfaces into
complex Grassman manifolds.
J. Diff. Geom. 27 (1988), 161-178.
\end{thebibliography}
\end{document}